\documentclass[12pt]{amsart}

\usepackage{epsf,amssymb}

\newtheorem{thm}{Theorem}[section]
\newtheorem{prop}[thm]{Proposition}
\newtheorem{defn}[thm]{Definition}
\newtheorem{lemma}[thm]{Lemma}
\newtheorem{cor}[thm]{Corollary}

\newtheorem{rmk}[thm]{Remark}
\newtheorem{claim}[thm]{Claim}

\newcommand{\R}{\mathbb{R}}
\newcommand{\Z}{\mathbb{Z}}

\newcommand{\N}{\mathbb{N}}
\renewcommand{\P}{\mathbb{P}}
\newcommand{\bdry}{\partial}
\newcommand{\s}{\vskip.1in}
\newcommand{\n}{\noindent}

\newcommand{\GCS}{\mathcal{GCS}}

\newcommand{\be}{\begin{enumerate}}
\newcommand{\ee}{\end{enumerate}}

\begin{document}
\title{Notes on the isotopy finiteness}

\author{Vincent Colin}
\address{Universit\'e de Nantes}
\email{Vincent.Colin@math.univ-nantes.fr}

\author{Emmanuel Giroux}

\address{\'Ecole Normale Sup\'erieure de Lyon}
\email{Emmanuel.GIROUX@umpa.ens-lyon.fr}

\author{Ko Honda}
\address{University of Southern California} 
\email{khonda@math.usc.edu}
\urladdr{http://math.usc.edu/\char126 khonda}

\date{This version: April 23, 2003.}


\begin{abstract}

{This is the less official, English version of the proof of the fact
  that every closed atoroidal $3$-manifold carries finitely many
  isotopy classes of tight contact structures.} 	

\end{abstract}
\maketitle

In this note we start where we left off in \cite{CGH1} and present the
proof of the following theorem: 

\begin{thm}\label{isotopy}
Let $V$ be a closed, oriented, atoroidal 3-manifold.  Then there exist only 
finitely many isotopy classes of tight contact structures on $V$. 
\end{thm}

Unless indicated otherwise, $V$ is a closed, oriented,
atoroidal 3-manifold, and contact structures on $V$ are {\em
cooriented}, i.e., oriented and positive. The space of tight contact
2-plane fields on $V$ is denoted by $Tight(V)$, and the set of isotopy
classes of tight contact structures on $V$ by $\pi_0(Tight(V))$.

\s\n
{\bf First Reduction:}  If $V$ is a connected sum $V_1\# V_2$, then a
gluing theorem of Colin \cite{Co1} implies that
$\pi_0(Tight(V_1))\times \pi_0(Tight(V_2))\stackrel \sim\rightarrow
\pi_0(Tight(V)).$  Since tight contact structures on $S^1\times S^2$
are also understood (there is a unique one up to isotopy), we may
assume that $V$ is irreducible.

\section{Weights on branched surfaces} 

\n In \cite{CGH1}, we constructed a finite number of pairs 
$(\mathcal{B}_i, \zeta_i)$, $i=1,\dots,k$, where:

\be
\item $\mathcal{B}_i$ is a branched surface, possibly with boundary,
\item $N(\mathcal{B}_i)\subset V$ is a branched surface neighborhood of 
$\mathcal{B}_i$,
\item $\zeta_i$ is a tight contact structure on $V\setminus N(\mathcal{B}_i)$,   
\ee
such that every $\xi\in Tight(V)$, up to isotopy, is {\em generated} by some   
$(\mathcal{B}_i,\zeta_i)$, i.e.,

\be
\item the fibers of $\pi_i: N(\mathcal{B}_i)\rightarrow \mathcal{B}_i$ are 
Legendrian,
\item $\zeta_i|_{V\setminus N(\mathcal{B}_i)}=\xi|_{V\setminus 
N(\mathcal{B}_i)}.$ 
\ee

\begin{rmk}  Note that a ``branched surface neighborhood''
  $N(\mathcal{B}_i)$ is usually not a neighborhood (in the topological
  sense) of $\mathcal{B}_i$, if we think of $\mathcal{B}_i$ as
  embedded in $V$.  Therefore, we will only think of
  $N(\mathcal{B}_i)$ as embedded inside $V$; $\mathcal{B}_i$ will be
  an abstract branched surface which is not embedded inside $V$,
  although there is a projection map $\pi_i:
  N(\mathcal{B}_i)\rightarrow \mathcal{B}_i$.
\end{rmk}    

\s\n
Fix some $(\mathcal{B}_i,\zeta_i)$ -- for simplicity, we omit the index $i$.
It suffices to prove the finiteness of isotopy classes of tight contact 
structures generated by $(\mathcal{B},\zeta)$.   Fix a tight contact structure 
$\xi_0$ generated by $(\mathcal{B},\zeta)$.  If $\xi$ is another tight contact 
structure generated by $(\mathcal{B},\zeta)$, then $\xi$ and $\xi_0$ differ only 
in the number of twists along each fiber of $\pi$ (which is Legendrian for both 
$\xi_0$ and $\xi$).  

\begin{lemma}
If $L$ is the branch locus of $\mathcal{B}$, then on each connected
component $B$ of $\mathcal{B}\setminus L$ the difference in the number
of twists is constant and is an integer.
\end{lemma}

\begin{proof}
Since $\xi$ and $\xi_0$ both agree on $\bdry N(\mathcal{B})$ and each fiber 
$\pi^{-1}(pt)$ is Legendrian for both, the difference in the number of twists is 
an integer.  Also, $\pi^{-1} B= B\times I$ is fibered by Legendrian arcs
$\{pt\}\times I$ and $B\times\{0,1\}$ is fixed for both $\xi_0$ and $\xi_1$.  
Continuity then guarantees that the integer does not vary over $B$. 
\end{proof}

\begin{figure} [ht]			
{\epsfysize=1.23in\centerline{\epsfbox{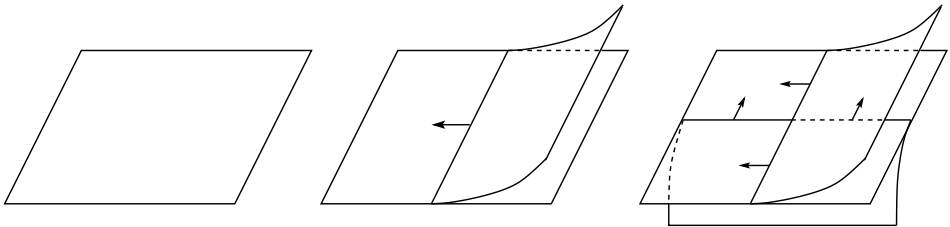}}} 		
\caption{The three local models for a branched surface.
The arrows indicate the branching direction.} 				
\label{branched}
\end{figure}

A connected component of $\mathcal{B}\setminus L$ will be called a {\em
sector}.  If $Q$ is the set of double points of $L$, then each
component $c$ of $L- Q$ is endowed with a normal direction in
$T\mathcal{B}$, called the {\em branching direction}, with the
property that for two of the sectors touching $c$, the branching
direction is the outward normal, and for the remaining sector touching
$c$, the branching direction is the inward normal. See Figure~\ref{branched}.

\begin{defn}
A {\em weight function} $w$ on $\mathcal{B}$ is a locally constant
function $w: \mathcal{B}\setminus L\rightarrow \Z$ (hence constant on
each sector) which satisfies the {\em branch equations}
$w(B_1)+w(B_2)=w(B_3)$, whenever $B_1$, $B_2$, $B_3$ are sectors
adjacent along a component 
$c$ of $L-Q$, and for $B_1$ and $B_2$ the branching direction is
outward and for $B_3$ the branching direction is inward.  We say $w$
is {\em positive} (resp.\ {\em nonnegative}) if
$w(\mathcal{B}\setminus L)\subset \Z^+$ (resp.\ $\subset \N$). 
\end{defn}

Therefore, to each $\xi$ generated by $(\mathcal{B},\zeta)$, we assign a weight 
function $w_\xi$ which on the connected component $B$ is the difference 
$t(\pi^{-1}(p), \xi_0) - t(\pi^{-1}(p),\xi)\in \Z$, where $t(\delta,*)$ is the 
``twisting number'' of the Legendrian arc $\delta$ with respect to the contact 
structure $*$, and  $p\in B$.

\begin{lemma}  \label{lemma: weight}
Let $w:\mathcal{B}\setminus L\rightarrow \Z$ be a nonnegative weight function.  
Then $w$ uniquely determines, up to isotopy rel $V\setminus N(\mathcal{B})$, a 
contact structure $\xi$ generated by $(\mathcal{B},\zeta)$. 
\end{lemma}

\begin{rmk} 
Of course it is difficult to tell which nonnegative weight 
functions $w$ correspond to tight contact structures. 
\end{rmk} 

\begin{proof}
The lemma is a consequence of the contractibility of $Diff^+(I)$, the
space of orientation-preserving diffeomorphisms of the unit interval.  

Let $\xi$ and $\xi'$ be two contact structures on $D^2\times I$ 
with coordinates $((x,y),t)$, which coincide on $D^2\times\{0,1\}$, have 
Legendrian fibers $\{pt\}\times I$, and have the same weight.  Then
they are given by $\alpha=\cos f_0(x,y,t)dx -\sin f_0(x,y,t)dy$ and
$\alpha'=\cos f_1(x,y,t)dx -\sin f_1(x,y,t) dy$ satisfying ${\bdry f_i\over
\bdry t}>0$ and $f_0(x,y,j)=f_1(x,y,j)$, where $i,j=0,1$.  By the
contractibility of $Diff^+(I)$, there
exists a 1-parameter family of functions $f_s:D^2\times I\rightarrow
\R$, $s\in[0,1]$, which 
satisfy ${\bdry f_s\over \bdry t}>0$ and are independent of $s$ on
$D^2\times\{0,1\}$.  Therefore $\xi$ and $\xi'$ are isotopic relative
to $D^2\times \{0,1\}$ through contact structures which have
$\{pt\}\times I$ as Legendrian fibers.  

In order to prove the lemma, we relativize the above discussion. Let
$\xi$ and $\xi'$ be two contact structures generated by $\mathcal{B}$
and which have the same weight.  If $B$ is a sector of $\mathcal{B}$,
then $\bdry B$ is a polygon $\delta_1\cup\delta_2\cup\dots\cup\delta_m$,
where $\delta_i$ are edges and the consecutive edges meet along triple
branch points of $\mathcal{B}$.  Since $(B\times I)\cap \bdry_v
N(\mathcal{B})$ consists of a union of $\delta_i\times
[a_i,b_i]$ where $[a_i,b_i]\subset [0,1]$, we are considering $\xi$ and $\xi'$ which
agree on $\delta_i\times [a_i,b_i]$.  We first line up $\xi$ and
$\xi'$ along $\delta_i\times ([0,1]-[a_i,b_i])$ using the
contractibility of $Diff^+(I)$, and then line up $\xi$ and $\xi'$ inside
$B\times I$, relative to $(\bdry B)\times I.$  This proves the lemma.
\end{proof}

Let $\mathcal{T}(\mathcal{B},\zeta)$ be the isotopy classes of tight contact
structures which are generated by $(\mathcal{B},\zeta)$.  (Recall that
$\mathcal{B}$ may have nonempty boundary.) The following is a
frequently used Amputation Lemma: 

\begin{lemma}[Amputation] 
If there exists a point $x\in \mathcal{B}$ such that the twisting number 
$\min_{\xi\in \mathcal{T}(\mathcal{B},\zeta)}
t(\pi^{-1}(x),\xi)>-\infty$, then we may ``amputate''  
all sectors of $\mathcal{B}$ containing $x$ from 
$\mathcal{B}$ to obtain a new branched surface $\mathcal{B}'$ (possibly with 
boundary) and tight contact structures $\zeta_1,\dots,\zeta_m$ such that 
$\mathcal{T}(\mathcal{B},\zeta)$ is generated by the $(\mathcal{B'},\zeta_i)$. 
\end{lemma}

\begin{proof}
Let $B$ be a sector of $\mathcal{B}$ which contains $x$.  
For any $y\in B$ and $\xi\in \mathcal{T}(\mathcal{B},\zeta)$, 
$t(\pi^{-1}(x),\xi)$ is also bounded below; hence there 
are only finitely many possibilities for $w_\xi(B)$.  For each $c\in \Z^+$, there 
is a finite number of tight contact structures $\zeta_i$ on
$V\setminus \pi^{-1}(\mathcal{B}-B)$ which agree with $\zeta$ on
$V\setminus N(\mathcal{B})$ and which admit Legendrian fibrations on
$\pi^{-1}(B)$ with $c=w_{\zeta_i}(B)$. Here, $\zeta_i$ is not unique
even if we fix $c$ because of the following: If  
$\delta$ is an edge of the polygon $\bdry B$ (described in Lemma~\ref{lemma: 
weight}), then $\bdry_v N(\mathcal{B})\cap (B\times I)$ may contain 
$\delta\times [a,b]$, where $[a,b]\subset [0,1]$.  Observe that there are finitely 
many ways of partitioning $$t(\{x\}\times [0,1])= t(\{x\}\times [0,a]) + 
t(\{x\}\times[a,b]) +t(\{x\}\times [b,1]),$$ where $t(\{x\}\times[a,b])$ is fixed. 
We can add an integer $d$ to $t(\{x\}\times [0,a])$ (resp.\  add $-d$ to 
$t(\{x\}\times [b,1])$), subject to the condition that $t(\{x\}\times [0,a])$ 
and $t(\{x\}\times [b,1])$ remain negative.
\end{proof}

By the following corollary, we may reduce to the case where
$\mathcal{B}$ has no boundary.

\begin{cor}
Suppose $\mathcal{T}(\mathcal{B},\zeta)$ is generated by
$(\mathcal{B},\zeta)$. If the branched surface $\mathcal{B}$ has
boundary, then $\mathcal{T}(\mathcal{B},\zeta)$ is also generated by a
finite collection $(\mathcal{B}_1,\zeta_1),\dots,
(\mathcal{B}_m,\zeta_m)$, where $\mathcal{B}_i$ are branched surfaces
without boundary. 
\end{cor}

\begin{proof}
Any point on $\bdry \mathcal{B}$ satisfies the conditions of the Amputation 
Lemma.  Also, each amputation reduces the complexity of $\mathcal{B}$, measured 
by the number of connected components of $\mathcal{B}\setminus L$.  Therefore,  
a finite number of applications of the Amputation Lemma yields the desired 
result.
\end{proof}

We now have the following simplification:

\begin{prop}   \label{prop: simplify}
There exist finitely many pairs $(\mathcal{B}_i,\zeta_i)$, $i=1,\dots,k$, such 
that:
\be
\item $\mathcal{B}_i$ is a branched surface {\em without} boundary.
\item $\zeta_i$ is a tight contact structure on $V\setminus N(\mathcal{B}_i)$.

\item Every $\xi_j\in \mathcal{T}(\mathcal{B}_i,\zeta_i)$, 
up to isotopy, is generated by $(\mathcal{B}_i, \zeta_i)$.  The corresponding 
weight function $w_{\xi_j}$ is sufficiently positive, i.e.,
$w_{\xi_j}(B)>>0$ for all sectors $B$ of $\mathcal{B}_i$.  
\item $\cup_{i=1}^k \mathcal{T}(\mathcal{B}_i,\zeta_i)=\pi_0(Tight(V))$.

\item  $\sup_j w_{\xi_j}(B)=+\infty$ for all sectors $B$ of
  $\mathcal{B}_i$. 

\ee
\end{prop}

\begin{rmk} 
It is possible that $\mathcal{B}_i$ is empty, in which case 
$\mathcal{T}(\mathcal{B}_i,\zeta_i)$ consists of just one element
$\zeta_i$. 
\end{rmk}

\begin{proof}
This follows from the Amputation Lemma.  Given a tight contact
structure $\xi$ generated by $(\mathcal{B},\zeta)$ for which there
exists a sector $B$ with insufficiently positive $w_{\xi}(B)$, we amputate $B$ to 
obtain $\mathcal{B}'$ (which we may assume has no boundary, after
further amputations) and transfer $\xi$ from $\mathcal{T}(\mathcal{B},\zeta)$ to 
$\mathcal{T}(\mathcal{B}',\zeta')$.  Moreover, if $\sup_j w_{\xi_j}(B)<+\infty$ 
for some $B$, then we can amputate $B$ from $\mathcal{B}$.
\end{proof}

\section{Surfaces carried by the branched surface $\mathcal{B}$}

Let $(\mathcal{B},\zeta)$ be one of the $(\mathcal{B}_i,\zeta_i)$ from 
Proposition~\ref{prop: simplify}; assume $\mathcal{B}$ is nonempty.  Each positive 
weight function $w_\xi$ corresponds to a closed surface $\mathcal{T}$ which is 
fully carried by $N(\mathcal{B})$.  Since $\mathcal{T}$ is 
transverse to the Legendrian fibration of $N(\mathcal{B})$, it follows that
$\mathcal{T}\pitchfork \xi$.   
This implies that each component of $\mathcal{T}$ is either a torus or a Klein 
bottle.  The following proposition allows us to restrict attention to the 
case where there are no Klein bottles.

\begin{prop}
Let $V$ be a closed, oriented, irreducible, and atoroidal 3-manifold which 
contains a Klein bottle $K$.  Then $|\pi_0(Tight(V))|<\infty$.
\end{prop}

\begin{proof}
Let $K$ be a Klein bottle, and let $\rho:N(K)\rightarrow K$ be its tubular 
neighborhood.  Then $\bdry N(K)$ is a torus, and since $\bdry N(K)$ is 
compressible, there exists a compressing disk $D$ for $\bdry N(K)$ in 
$V\setminus N(K)$.  Now, irreducibility implies that either $N(D)\cup N(K)$ is a 
3-ball or $V\setminus N(K)$ is a solid torus.  The first option is not possible, 
since $K$ does not separate $B^3$, whereas every closed surface in $B^3$ must 
separate.  The latter implies that $V= N(K)\cup (S^1\times D^2)$.
Now, $N(K)$ admits a Seifert fibration over the disk with two singular
fibers, where the regular fibers are isotopic to the boundary of a
M\"obius band ($\rho^{-1}(\delta)$ for an appropriate $\delta\subset
K$).  If the meridional curve of $S^1\times D^2$ is a
regular fiber of $N(K)$, then $V$ is $\R\P^3\# \R\P^3$, which has a
unique tight contact structure.  Otherwise,
$V$ is a Seifert fibered space over $S^2$ with three singular fibers
and Seifert invariants $({1\over 2}, {1\over 2}, {\beta\over
\alpha})$.  If ${\beta\over \alpha}$ is an integer, $V$ is a lens
space and it is well-known that $|\pi_0(Tight(V))|< \infty$.
Otherwise, the finiteness of tight contact structures on these {\em
  prism manifolds} follows from Section~\ref{small Seifert}. 
\end{proof}


From now on we assume that every connected surface carried by $N(\mathcal{B})$
is a torus.

\begin{lemma} There exists a surface $\mathcal{T}$ fully carried by
  $N(\mathcal{B})$ such that $\mathcal{T}\supset \bdry_h
  N(\mathcal{B})$, and there is an orientation on $\mathcal{T}$ which
  agrees with the normal orientation on $\bdry N(\mathcal{B})$ along
  $\mathcal{T}\cap \bdry_h N(\mathcal{B})$.
\end{lemma}

\begin{proof}
By doubling $\mathcal{T}$ if necessary, we ensure that each fiber of
$N(\mathcal{B})$ intersects $\mathcal{T}$ at least twice.  If
$\Sigma\subset \bdry_h N(\mathcal{B})$ is a connected component, then
we can flow $\Sigma$ along the Legendrian fibers until we hit
$\mathcal{T}$ for the first time. Reversing this process, we can pull
that portion of $\mathcal{T}$ to the boundary component $\Sigma$. 

If, on any component $T$ of $\mathcal{T}$, there are at
least two components of $T\cap \bdry_h N(\mathcal{B})$ and the
orientations on $T\cap \bdry_h N(\mathcal{B})$ are inconsistent, then
we double $T$ to obtain tori $T_+$ and $T_-$, and place components of
$\bdry_h N(\mathcal{B})$ of one orientation onto $T_+$ and components
of the opposite orientation onto $T_-$.
\end{proof}

\section{The Degree Lemma}

Consider the pair $(\mathcal{B}, \zeta)$.
Let $A$ be a connected component of $\bdry_v N(\mathcal{B})$ and
$c$ be a boundary component of $A$.  Define $\deg(A)$, the {\em 
degree} of $A$ to be the absolute value of the degree of the image
of $\zeta_x$ with respect to $T_xA$  in the quotient
$T_xV/T_x(\mbox{Legendrian fiber})$, as 
$x$ ranges over $c$. Here, we use the absolute value because there is
not necessarily a coherent orientation on the fibers of $\pi:
N(\mathcal{B})\rightarrow \mathcal{B}$.

\begin{claim}\label{degree}
We may assume that for each component $A$ of $\bdry_v N(\mathcal{B})$
and each component $c$ of $\bdry A$, there are exactly $2\deg(A)$
points along $c$ where $T_x A= \zeta_x$.  
\end{claim}

\begin{proof}
Extend $A$ along the Legendrian fibers to an annulus $A'$, where
$A'-A\subset N(\mathcal{B})$, $\xi|_{A'}$ agree for all $\xi\in
\mathcal{T}(\mathcal{B},\zeta)$, and the condition of the claim holds
for $\bdry A'$. (This is made possible by all $\xi\in
\mathcal{T}(\mathcal{B},\zeta)$ having sufficiently large $w_{\xi}$).

The claim follows from ``making a slit'' along $A-A'$. More precisely,
consider a thickened annulus $A'\times I \subset N(\mathcal{B})$, where
$A'\times\{1\}=A'$ and each $\xi$ is $I$-invariant on the thickened
annulus. (In particular, each $A'\times\{t\}$ is fibered by Legendrian
intervals.)  Then the new $N(\mathcal{B})$ is $N(\mathcal{B})\setminus
(A'\times I)$ with some corners rounded, so that $\bdry_v
N(\mathcal{B})=A'\times\{0\}$. 
\end{proof}

From now on, all components of $\bdry_v N(\mathcal{B})$ are assumed to
satisfy the above claim.

\begin{lemma}[Degree Lemma]
After possible amputations, we may assume that every connected
component $A$ of $\bdry_v N(\mathcal{B})$ with nonzero degree
intersects $\mathcal{T}$ along homotopically trivial curves on
$\mathcal{T}$.
\end{lemma}

\begin{proof}
Let $T$ be a component of $\mathcal{T}$ which intersects a component
$A$ of $\bdry_v N(\mathcal{B})$ with $\deg(A)\not=0$ along $c$. Suppose
$c$ is homotopically essential on $T$. For any $\xi\in
\mathcal{T}(\mathcal{B},\zeta)$ with $w_\xi$ 
sufficiently positive, there exists an embedding
$\phi:T\times[0,1]\rightarrow N(\mathcal{B})$, where $\phi(T,0)=T$
and
$\phi^* \xi$ is given by $\cos (f(x,y)+2\pi t)dx -\sin (f(x,y)+2\pi t)
dy=0$. Here, the coordinates on
$T\times[0,1]=\R^2/\Z^2\times[0,1]$ are $(x,y,t)$, and $f$ is a
circle-valued function $T\rightarrow \R/2\pi\Z$.  We then have the following:

\begin{claim}\label{perturb}
Inside $\phi(T\times[0,1])$ there exists a torus
$T'$ isotopic to $T$ and transverse to the Legendrian fibers,
such that $T'$ is convex and
$\#\Gamma_{T'}\leq 2 \deg(A)$. 
\end{claim}

\begin{proof}[Proof of Claim~\ref{perturb}]
After a $C^\infty$-small perturbation of $T$, we may assume $T$ is
convex.  Since $T\pitchfork \xi$ for any $\xi\in
\mathcal{T}(\mathcal{B},\zeta)$, the characteristic foliation $\xi T$ is
nonsingular, and hence $\#\Gamma_T$ is the same as the number of
closed orbits $\gamma_i$ of $\xi T$.  $T\setminus \cup_i \gamma_i$ are
annular components which are either Reeb (no transverse arc with
endpoints on the boundary which intersects every leaf) or taut (there
exists such a transverse arc). We may assume $c$ is transverse to
$\cup_i\gamma_i$.  By inspecting the connected components of
$c\setminus \cup_i\gamma_i$, every (separating or nonseparating) arc
inside a Reeb component contributes at least one tangency, whereas arcs
inside taut components do not necessarily contribute.  Therefore, the
number of Reeb components is bounded above by $2\deg (A)$ ($=$ the
number of tangencies of $c$), if the $\gamma_i$ have nontrivial
geometric intersection with $c$.  To see that the $\gamma_i$ have
nontrivial geometric intersection with $c$, observe that
$2\deg(A)$, the signed count of tangencies of $c$ and $\xi T$, is
invariant under isotopy. If the geometric intersection number is zero,
then the degree must be zero. Finally,  all the taut components can be
removed by isotoping $T$ a bounded distance within $\phi(T\times[0,1])$. 
\end{proof}

The key feature of the convex torus $T$ modified as in the above lemma
is that $\#\Gamma_T$ is bounded independently of the choice of $\xi
\in\mathcal{T}(\mathcal{B},\zeta)$.  Suppose that $\xi\in \mathcal{T}
(\mathcal{B},\zeta)$ satisfies $w_\xi >> nw_T$, where $n=\deg
(A)={1\over 2} \#\Gamma_T$. Then there exists an embedding $\psi:
T\times[0,n]\rightarrow N(\mathcal{B})$, where $T\times\{0\}=T'$ and
$\psi^*\xi$ is given by $\cos (g(x,y)+2\pi t)dx -\sin (g(x,y) +2\pi t)
dy=0$.  If we could remove $\psi(T\times[0,n])$ and reglue
$\psi(T\times\{0\})$ with $\psi(T\times\{n\})$ via the natural
identification given by the Legendrian fibration, we obtain a contact
structure $\xi'$ corresponding to the weight $w_\xi - nw_T$.  Now,
$\xi$ and $\xi'$ are isomorphic, since they differ by Dehn twists
along tori.  They are isotopic, since either (i) $T$ bounds a solid
torus or (ii) $T$ bounds a knot complement inside $B^3$, and we can
use the fact that a diffeomorphism of $B^3$ relative to the boundary
is isotopic to the identity rel boundary.  Therefore, we may
inductively reduce $w_\xi \rightarrow w_\xi - nw_T$ until some sector
$B$ has small weight. Such a sector $B$ can be amputated.
\end{proof}

With the Degree Lemma at hand, we prove the following useful proposition:

\begin{prop}\label{improved degree lemma}
After possible amputations, we may assume that, given a connected
component $A$ of $\bdry_v N(\mathcal{B})$, 
\be
\item $\deg (A)=0$ if and only if both components of $\bdry A$ are
  essential on $\mathcal{T}$.
\item $\deg (A)=1$ if and only if both components of $\bdry A$
  bound disks in $\mathcal{T}$.
\ee
\end{prop}

\begin{proof}
By the Degree Lemma, if a component of $\bdry A$ is essential, then
$\deg(A)=0$.  On the other hand, if a component $c$ of $\bdry A$ bounds a
disk $D$ in $\mathcal{T}$, then, by Claim~\ref{degree} and
the nonsingularity of the characteristic foliation on $D$, there
can only be two points along $c$ where $T_xA=\zeta_x$.  Hence $\deg
(A)=1$.  Therefore, either both components of $\bdry A$ are essential,
or both bound disks.  The proposition follows.
\end{proof}

\begin{rmk}
Observe that if both components of $\bdry A$ bound disks, then the
disks must both be on the same side of $A$.
\end{rmk}

\section{Elimination of disks of contact}

In this section, we simplify the branched surface neighborhood
$N(\mathcal{B})$ by eliminating disks of contact.  A {\em disk of
  contact} is a properly embedded disk $D\subset N(\mathcal{B})$
transverse to the fibers of $N(\mathcal{B})$,
whose boundary is on $\bdry_v N(\mathcal{B})$.

\begin{lemma}\label{replace}
Let $A$ be a component of $\bdry_v N(\mathcal{B})$. If there exists a
disk of contact $D$ with boundary on $A$, then the boundary components
$c_1$, $c_2$ of 
$A$ bound disks $D_1$, $D_2\subset \mathcal{T}$ so that $D_i$ is in the
interior of $N(\mathcal{B})$ near $\bdry D_i$.
\end{lemma}

\begin{proof}
Since $A$ admits a disk of contact $D$ and the characteristic
foliation on $D$ is nonsingular, $\deg(A)$ must equal one.  
By the Degree Lemma,
$c_i$ must bound a disk $D_i$ in $\mathcal{T}$. Note that $D_i$ cannot
be in the ``opposite direction'' from $D$, namely $D_i$ cannot contain
the component of $\bdry_h N(\mathcal{B})$ adjacent to $c_i$.  Otherwise
$D\cup D_i$ (together with some pieces of $A$, and after some
rounding) will form an immersed 2-sphere transverse to the Legendrian
fibration, a contradiction. 
\end{proof}

\begin{rmk}
It is conceivable that $D_1\subset D_2$ or vice versa.
\end{rmk}

\begin{prop}\label{no disks of contact}
Let $V$ be a closed, atoroidal, irreducible manifold. There exists a
finite number of pairs $(N_i,\zeta_i)$, $i=1,\dots,k$, satisfying the
following: 
  
\be
\item $N_i\subset V$ is a finite union of thickened tori $T^2\times
  [0,1]$ and annuli $A\times [0,1]$, where each $A\times\{j\}$,
  $j=0,1$, is glued {\em essentially} onto some $\bdry(T^2\times
  [0,1])$, and some boundary components of the $T^2\times I$ may be identified.  

\item $\zeta_i$ is a tight contact structure on $V\setminus N_i$.

\item If $\mathcal{T}(N_i,\zeta_i)$ is the set of isotopy classes of
  tight contact structures $\xi$ on $V$ which agree with $\zeta_i$ on
  $V\setminus N_i$, then $\cup_i
  \mathcal{T}(N_i,\zeta_i)=\pi_0(Tight(V))$.
\ee 
\end{prop}

\begin{proof}
We first eliminate all disks of contact from $N(\mathcal{B})$, while
preserving the condition that $N(\mathcal{B})$ fully carry a union of
tori $\mathcal{T}$. (See Remark~\ref{loss} below.) If there is a disk
of contact for $A$, then using Lemma~\ref{replace}, we may replace it
with disks $D_1$ and $D_2$ of contact (also for $A$) in
$\mathcal{T}$. Without loss of generality, assume $D_1$ is an
innermost disk of contact for $\mathcal{T}$. Then either $D_1$ and
$D_2$ are disjoint, or $D_1\subset D_2$.  Since $D_1$ may contain some
disks of $\bdry_h N(\mathcal{B})$, we make a small isotopy of $D_1$
along the fibers, to push $D_1$ away from $\bdry_h N(\mathcal{B})\cap 
int(D_1)$.  Call the isotoped disk $D_1'$. Then modify
$N(\mathcal{B})\rightarrow N(\mathcal{B})\setminus D_1'$,
$\mathcal{T}\rightarrow (\mathcal{T}\setminus D_1)\cup D_1'$, and
$D_2\rightarrow (D_2\setminus D_1)\cup D_1'$ if $D_1\subset D_2$
(rename them $N(\mathcal{B})$, $\mathcal{T}$, and $D_2$). 

We will now explain how to realize $D_1'$ as a convex surface (with
Legendrian boundary), so that the tight contact structure $\zeta$ on
$N(\mathcal{B})$ extends uniquely to a contact structure on 
$N(\mathcal{B})\cup N(D_1')$.  After rounding the corners of $\bdry
N(\mathcal{B})$ and perturbing, $\bdry N(\mathcal{B})$ becomes
convex.  The fact that $\deg(A)=1$ translates into $tb(\bdry
D_1')=-1$, when we realize $\bdry D_1'$ as a Legendrian curve on
$\bdry N(\mathcal{B})$.  Now, if $D_1'$ is perturbed into a convex
surface with Legendrian boundary, there is only one possibility (up to
isotopy) for $\Gamma_{D_1'}$.  Hence, after applying the Flexibility
Theorem, we may assume that for any $\xi\in
\mathcal{T}(N(\mathcal{B}),\zeta)$, $D_1'$ can be taken to have the
same characteristic foliation.

Since (the new) $\mathcal{T}$ is not fully carried by (the new)
$N(\mathcal{B})$, we modify $\mathcal{T}$ as follows: Let $T$ be the
connected component of $\mathcal{T}$ containing $D_2$, and let $T'$ be
a parallel push-off.  Surger $T\rightarrow (T\setminus D_2)\cup
A\cup D_1'$ and round.  Doubling the surgered torus, we obtain a fully
carried $\mathcal{T}$ containing $\bdry_h N(\mathcal{B})$. Since there
are only finitely many components of $\bdry_h N(\mathcal{B})$, we can
eliminate all the disks of contact.  Observe that components of
$\bdry_v N(\mathcal{B}) \cap \mathcal{T}$ which were homotopically
essential (resp.\ homotopically trivial) remain essential (resp.\
trivial) after the surgery and doubling operations.

Having eliminated all disks of contact, we now examine the connected
components of $\bdry_h N(\mathcal{B})$.  Indeed, there are only three
possibilities: (i) disks, (ii) annuli which are essential on
$\mathcal{T}$, and (iii) tori.  All the disk components of $\bdry_h
N(\mathcal{B})$ can be eliminated as follows: Let $D$ be a disk component of $\bdry_h N(\mathcal{B})$ and $A$ an
annulus of $\bdry_v N(\mathcal{B})$ which shares a boundary component
with $D$. Then $\deg(A)$ must be nonzero, and if $S$ is a component of
$\bdry_h N(\mathcal{B})$ which intersects the other boundary component
of $A$, then, by the Degree Lemma, $S$ cannot be a homotopically
essential annulus. Hence, $S$ is also a disk.  Now, $D\cup A\cup S$ is
a sphere which bounds a 3-ball $B^3$ on one side or another.  In one
case, we take $N(\mathcal{B})\cup B^3$, and in the other case we take
$N(\mathcal{B})\setminus B^3$.  Eventually, the horizontal disk
components are removed.  This implies that all the components of
$N(\mathcal{B})\setminus \mathcal{T}$ are thickened tori or thickened
annuli which are glued essentially onto the boundary of the thickened tori.
\end{proof}

\begin{rmk}\label{loss}
In eliminating disks of contact, we lose control over the Legendrian
fibration, although the topological fibration still exists.
Therefore, instead of isotopy classes of tight contact structures
which are {\em generated} by a pair $(N(\mathcal{B}),\zeta)$, we must
consider isotopy classes of tight contact structures on $V$ which
simply agree with $\zeta$ on $V\setminus N$.  Due to this loss of
information, we must repeat the finiteness study for simpler $N$ and
$V$, namely when $V$ is a small Seifert space.  This study will be
conducted in the next two sections.
\end{rmk}

\section{Reduction to the small Seifert case}

Let $(N,\zeta)=(N_i,\zeta_i)$ be a pair as in Proposition~\ref{no disks
  of contact}.

\begin{lemma}
$N$ is a graph manifold with nonempty boundary.
\end{lemma}

\begin{rmk}
Our graph manifolds may be disconnected, and the components may be
Seifert fibered spaces.
\end{rmk}

\begin{proof}
Suppose $N=V$.  Then $N$ consists only of $T^2\times I$ components,
glued successively to give a torus bundle over $S^1$.  Therefore $V$
is toroidal, a contradiction.  

Now, whenever two $T^2\times I$ components share a common boundary (they
cannot share both boundary components), they can be merged into a
single $T^2\times I$. Next, if we cut $N$ along the union of
$T^2\times \{{1\over 2}\}$, then the connected components are
diffeomorphic to $S^1$ times a compact surface with boundary.
Therefore $N$ is a graph manifold. 
\end{proof}

\begin{prop}
If $\mathcal{T}(N,\zeta)$ is not finite, then $V$ is a Seifert fibered
space over $S^2$ with 3 singular fibers.
\end{prop}

\begin{proof}
Let $T$ be a boundary component of $N$, which is a graph manifold.
Since $T$ is compressible, there is a compressing disk $D$ for $T$.  By using
an innermost argument and switching $T$ if necessary, we may assume
that $D\subset N$ or $D\subset V\setminus N$. 

Suppose $D\subset N$.  Since $N$ is a graph manifold, it is
irreducible.  Hence $T$ bounds a solid torus $W$ which contains
components of $N$.  By the finiteness of tight contact structures on $W$, all
components of $N$ inside $W$ can be removed from $N$ without affecting
the infiniteness of $\mathcal{T}(N,\zeta)$.

Now suppose $D\subset V\setminus N$.  By the irreducibility of $V$,
either $T$ bounds a solid torus $W$
in $V\setminus N$, or $T\cup D$ is contained in a 3-ball whose
boundary lies outside $N$.  In the latter case, we may throw away all
the components of $N$ inside the 3-ball (including the
component of $N$ which is bounded by $T$).  In the former case, we
consider $N\cup W$. Let $M$ be the maximal Seifert fibered component
of $N$ with $M\cap W\not=\emptyset$, as given by the canonical torus
(Jaco-Shalen-Johannson) decomposition.  Also let $\pi: M\rightarrow S$ be
the projection onto $S$, a compact surface with boundary.  If $S$ is a
disk with at most one singular point, then $M$ is a solid torus and
$V$ is a lens space.  If $S$ is an annulus without any singular
points, then $M=T^2\times I$ is a connected component of $N$, and
$M\cup W$ is a solid torus with a finite number of tight contact
structures, hence can be excised from $N$.  

For any other $S$, if the meridian of $W$ does not bound a regular
fiber in $M$, then $M\cup W$ is a maximal Seifert fibered component of
$N\cup W$.  In this case, $N\cup W$ is a graph manifold -- we will
rename this $N$.  If the meridian of $W$ does bound, then let $c$ be a
boundary component of $S$ which 
corresponds to $T$.  There exists an
arc $d\subset S-\{\mbox{singular points}\}$ with endpoints on $c$,
which is not $\bdry$-parallel in $S-\{\mbox{singular points}\}$. Now,
the union of $\pi^{-1}(d)$ and two meridional disks of $W$ is a
2-sphere $K$, which must bound a 3-ball $B^3$ on one side, by irreducibility.
This implies, first of all, that $S$ is a planar surface; otherwise
there exists an arc $d$ and a closed curve $\delta$ on $S$ which
intersect precisely once, contradicting the fact that $K$
separates.  Next, consider $M'=M\cup W\cup B^3$.  One of the components
of $S\setminus d$ is contained in $B^3$, and the other is a planar
surface $S'$ with fewer boundary components.  Since $K$ bounds a
3-ball, $S'$ bounds a solid torus in $M'$.  Renaming $N\cup W\cup B^3$
and $M\cup W\cup B^3$ to be the new $N$ and $M$, and continuing in
this manner, we see that $M\cup W$ is a solid torus (or $V$ is a lens space).

By repeating the above argument, we inductively reduce the number of
connected components of $\bdry N$, while ensuring that $N$ is either a graph
manifold or the empty set (which is disallowed by hypothesis).  Since
$V$ is atoroidal, $V=N$ would then be a small Seifert fibered space or
a lens space (which is also disallowed).
\end{proof}

\section{The small Seifert case}\label{small Seifert}

Let $V$ be a small Seifert space, i.e., a Seifert fibered space with 3
singular fibers over $S^2$.  The tubular neighborhoods of the singular
fibers $F_i$, $i=1,2,3$, are denoted by $V_i$.

\subsection{Case 1}

We restrict attention to the set of tight contact structures for which
there exists a Legendrian regular fiber with twisting number $=0$,
where the twisting is measured using the projection.

\begin{rmk} (Well-definition of twisting number) Take a neighborhood
  of a regular fiber, and on the boundary we look at curves which
  bound on either side.  If they are different, the framing is
  well-defined; if they are the same, then that means that there is a
  horizontal surface (after making the surface incompressible).  This
  works in our case.  You can also use the uniqueness of the Seifert
  fibration.
\end{rmk}

\begin{claim} Given a tight contact structure $\xi$ with a
  zero-twisting Legendrian regular fiber, there exists an
  isotopy of $\xi$ for which $\bdry V_i$ are convex, 
  $\Gamma_{\bdry V_i}$ are vertical, and $\#\Gamma_{\bdry V_i}=2$.
\end{claim} 

Let $W=V\setminus \cup_{i=1}^3 V_i$.  Let $\pi: W\rightarrow S$ be the
projection induced by the fibration, where $S$ is the 3-punctured
sphere.

\begin{claim} Let $\widehat S$ be a convex surface which is the image
  of a section $s:S\rightarrow W$.  Then $\Gamma_{\widehat S}$ consists
  of three nonseparating arcs, each of which connects distinct
  boundary components of $\widehat S$. 
\end{claim}

\begin{proof}
Suppose there is a separating arc $c$ on $\widehat S$. Then it is
$\bdry$-parallel.  Let $V_i$ be the solid torus to which the bypass
corresponding to $c$ is attached.  Then the thickening $V'_i$ of $V_i$
has convex boundary with $\Gamma_{\bdry V'_i}$ parallel
to $\bdry \widehat S\cap \bdry V_i$.  Now $V'_i$ can be thickened
again to $V''_i$ so as to contain Legendrian fibers with twisting
number zero (taken, for instance, by pushing off a zero twisting curve
on $V_j$, $j\not=i$).  Some intermediate convex torus between $\bdry
V'_i$ and $\bdry V''_i$ will then have dividing curves parallel to the
meridian of $V''_i$, and hence the contact structure on $V''_i$ is overtwisted.

Now, if $c$ is a closed component of $\Gamma_{\widehat S}$, it is
parallel to some boundary component $\bdry V_i \cap \widehat S$, and
the same argument as above shows that the contact structure is
overtwisted.
\end{proof}

Write $W=S\times S^1$.  If $\bdry S= c_1\sqcup c_2\sqcup c_3$, then
$\bdry W= \sqcup_{i=1}^3 c_i\times S^1$.

\begin{prop} There is a 1-1 correspondence between isotopy
classes of tight contact structures on $W=S\times S^1$ with fixed convex boundary
where $\Gamma_{c_i\times S^1}$, $i=1,2,3$, consists of $2k_i$ parallel
curves isotopic to the regular fiber, and isotopy classes of
multicurves on $S$ which have no homotopically trivial components and
which have $2k_i$ fixed endpoints on $c_i$.
\end{prop}

The proposition can be found in \cite{Gi5,Ho2}.

Observe that there are infinitely many isotopy classes of possible
dividing sets on $S$ relative to the boundary; however, non-relative
isotopy classes are finite in number (in fact there are
two). Moreover, if $\Gamma$ and $\Gamma'$ are two allowable dividing
sets on $S$ which are isotopic but not isotopic relative to the
boundary, they differ by Dehn twists parallel to the boundary components.
In other words, we may assume that $\Gamma=\Gamma'$ when restricted to $\widehat
S'=\widehat S-\cup_i A_i$, where $A_i$ is a collared neighborhood of
$\widehat S\cap \bdry V_i$. Let $W'=\pi^{-1}(\pi (\widehat S'))$, and
$V'_i$ be the component of $V\setminus W'$ containing $V_i$. 
 
Given a tight contact structure $\xi$ on $V$ with a zero-twisting
Legendrian regular fiber, there exists an isotopy of $\xi$ so that
$\xi|_{W'}$ is one of two types.   Now since there are only finitely
many isotopy classes of tight contact structures on solid tori with a
fixed boundary condition, we conclude that there are finitely many
isotopy classes of tight contact structures on $V$ with a Legendrian
regular fiber with zero twisting.


\subsection{Case 2}

In this case we only consider the set of tight contact structures on $V$
which do not have Legendrian regular fibers with zero twisting
number. 

\begin{rmk}  The Seifert fibered space with 3 singular fibers over
  base $S^2$ is denoted $({\beta_1\over \alpha_1}, {\beta_2\over
  \alpha_2}, {\beta_3\over \alpha_3})$.  ${\beta_i\over \alpha_i}$
  is then the slope of the meridional disk of $V_i$, seen from
  $W=S\times S^1$.  Here, $\bdry S\cap V_i$ has slope zero and the
  regular fibers have slope $\infty$.
\end{rmk}

We define ${\beta_i'\over \alpha_i'}$ with $GCD(\beta_i',\alpha_i')=1$ to
be the greatest rational number such that $\beta_i' \alpha_i-\beta_i
\alpha_i'=1$.  (Viewed on the Farey tessellation, ${\beta_i'\over
\alpha_i'}$ is the point closest to $+\infty$ on $({\beta_i\over
\alpha_i},+\infty)$ with an edge to ${\beta_i\over \alpha_i}$.) 

\begin{claim}\label{twist number}
  Suppose $\xi$ is a contact structure so that $V_i$ is the
  standard neighborhood of a Legendrian singular fiber and
  $\Gamma_{\bdry V_i}$ has slope in $({\beta_i\over \alpha_i},
  {\beta_i'\over \alpha_i'})$.  If there exists a bypass along a
  Legendrian regular fiber on $\bdry V_i$, then $\xi$ is isotopic to a
  contact structure $\xi'$ for which $V_i$ is the standard
  neighborhood of a Legendrian singular fiber with one higher twisting
  number.  
\end{claim}

Define $\GCS$ to be the isotopy classes of tight contact structures
for which there exists a representative $\xi$ satisfying the following
conditions:

\begin{enumerate}

\item $V_1$ and $V_2$ are standard neighborhoods of singular
  Legendrian fibers.

\item The annulus $A$ connecting $\bdry V_1$ to $\bdry V_2$ is convex,
  contains no $\bdry$-parallel dividing curves, and is fibered by
  Legendrian regular fibers with maximal twisting number.

\end{enumerate}


\begin{claim} There are only finitely many isotopy classes of tight
  contact structures which are not in $\GCS$.
\end{claim}

\begin{proof}
Let $F$ be a Legendrian regular fiber with maximal twisting number,
$V_i$, $i=1,2$, be standard neighborhoods of Legendrian singular
fibers with boundary slopes in $({\beta_i\over \alpha_i},
{\beta_i'\over \alpha_i'})$, and $A$ be a convex annulus from $\bdry
V_1$ to $\bdry V_2$ which contains $F$.  If $A$ has no
$\bdry$-parallel arc, then the contact structure $\xi$ is in $\GCS$.
Otherwise, we attach the corresponding bypass to thicken $V_i$ by
using Claim~\ref{twist number}.  Thus, we are reduced to considering the
case when the boundary slopes of $V_i$ are ${\beta_i'\over
\alpha_i'}$.

Consider $\gamma_i$, the shortest increasing path in the Farey
tessellation from ${\beta_i'\over \alpha_i'}$ to $+\infty$.  If the
slope of $\bdry V_i$ is in (the vertices of) $\gamma_i$, then
attaching a bypass corresponding to 
a $\bdry$-parallel component of $\Gamma_A$ produces a solid torus with
boundary slope which is the next term in the path $\gamma_i$.  Hence, repeating this
operation, we obtain a contact structure for which $A$ has no
$\bdry$-parallel dividing curves and the boundary slopes of $V_1$ and
$V_2$ are in the sequences $\gamma_1$, $\gamma_2$.

Since the vertices of $\gamma_i$ are finite in number, and
$\Gamma_{\bdry V_3}$ is determined by the above data, the proof
follows from
using the finiteness of tight contact structures on solid tori with
fixed boundary slopes and a fixed number of dividing curves.  
\end{proof}

We are now left to consider $\GCS$.

\begin{prop}
If $\sum_{i=1}^3 {\beta_i\over \alpha_i}\not=0$, then $\GCS$ is finite.
\end{prop}

\begin{proof}
We argue by contradiction.  Suppose there exists an infinite sequence
$\xi_k$ of tight contact structures in $\GCS$.  Then the boundary
slopes on $V_1$ and $V_2$ for $\xi_k$ converge to meridional slopes
${\beta_1\over \alpha_1}$ and ${\beta_2\over \alpha_2}$.  (In fact, if
the two boundary slopes remain bounded away from the meridional
slopes, finiteness follows from the finiteness of tight contact
structures on solid tori, and if one boundary slope tends to its
meridional slope, the other must also tend to its meridional slope
because of the connecting annulus $A$.)

Letting $V_3$ be the torus obtained from $V\setminus (V_1\cup V_2\cup A)$
by rounding edges, the boundary slopes $s_k$ of $\xi_k$ on $V_3$ tend
to $s=-{\beta_1\over \alpha_1}-{\beta_2\over \alpha_2}$.  Let $\gamma$
be the interval $({\beta_3\over \alpha_3},s)$, where if $s<
{\beta_3\over \alpha_3}$ we understand it to be $({\beta_3\over
\alpha_3}, +\infty] \cup [-\infty,s)$.  Now, if $s< {\beta_3\over
\alpha_3}$, then for $k$ large enough there exists a convex
torus $T$ in $V_3$ parallel to $\bdry V_3$ with slope $\infty$.
This is ruled out by the assumption of Case 2.  Therefore we
assume $s> {\beta_3\over \alpha_3}$.

Now consider $\xi_k$ where $s_k>s$.  Then there exists a convex torus
$T$ in $V_3$ parallel to $\bdry V_3$ with slope $s$.  For $k$ large
enough, $\Gamma_T$ intersects the regular fiber fewer times than
$\Gamma_{\bdry V_3}$; this contradicts the maximality of the twisting
number of the Legendrian regular fiber in $A$.

Finally, consider $\xi_k$ where $s_k<s$.  On the interval $\gamma$,
let $s'$ be the vertex of the Farey tessellation closest to ${\beta_3\over
\alpha_3}$ with an edge to $s$.   For $k$ large enough, $s_k$ is in the
interval $(s',s)$, and hence there exists a convex torus $T$ in $V_3$
with slope $s'$.  Its dividing curves intersect the regular fiber in
fewer points than those of $\bdry V_3$, which is again a
contradiction.
\end{proof}

Suppose now that $\sum_i {\beta_i\over \alpha_i} =0$.  

\s\n
{\em Convention:} We will normalize the Seifert invariants so that
$0<{\beta_1\over \alpha_1}, {\beta_2\over \alpha_2}<1$ and $-2<
{\beta_3\over \alpha_3} <0$, and we take the $\alpha_i$ to be positive
integers.

\begin{thm}\label{infinite GCS} 
  If $\GCS$ is infinite, then $V$ is an elliptic torus
  bundle over the circle, and hence is toroidal.
\end{thm}

\begin{rmk} Seifert fibered spaces over $S^2$ with three singular
  fibers which are torus bundles over the circle are classified --
  they have Seifert invariants $\pm (-{1\over 2}, {1\over 4}, {1\over
  4})$, $\pm (-{1\over 2}, {1\over 3}, {1\over 6})$, and $\pm
  (-{2\over 3}, {1\over 3}, {1\over 3})$.  They satisfy the property
  that $\sum_i{1\over \alpha_i}=1$.
\end{rmk}

Suppose $\GCS$ is infinite.  Then there exists an infinite number of
positive pairs $(k_1,k_2)$ so that
$k_1\alpha_1+\alpha_1'=k_2\alpha_2+\alpha_2'$, and tight contact
structures whose corresponding boundary slopes of $V_i$ are equal to
${k_i\beta_i+\beta_i'\over k_i\alpha_i+\alpha_i'}$.  Indeed, only
finitely many contact structures induce the given boundary slopes on
$V_1$ and $V_2$.  So the existence of infinitely many contact
structures in $\GCS$ implies the existence of infinitely many
$(k_1,k_2)$.

The solutions of the equation
$k_1\alpha_1+\alpha_1'=k_2\alpha_2+\alpha_2'$ are parametrized by a
number $k\in {1\over GCD(\alpha_1,\alpha_2)} \N$ in the following
way: Given a particular solution $(r_1,r_2)$, other solutions are
parametrized by $k_1=k\alpha_2+r_1$ and $k_2=k\alpha_1+r_2$.  Observe
that there exists a subsequence $k\rightarrow +\infty$ and tight contact
structures $\xi_k$ in $\GCS$.

We compute the boundary slope of $V_3$ to be:
\begin{eqnarray*}
s_k  & = & {1- ((k\alpha_2+r_1)\beta_1 +\beta_1')-
  ((k\alpha_1+r_2)\beta_2 +\beta_2')\over
  (k\alpha_2+r_1)\alpha_1+\alpha_1'}\\
& = & {k(-\alpha_2\beta_1-\alpha_1\beta_2) + (1-r_1\beta_1 -\beta_1'
    -r_2\beta_2 -\beta_2')\over k(\alpha_1\alpha_2) + r_1\alpha_1
    +\alpha_1'}
\end{eqnarray*}
where $(k\alpha_2+r_1)\alpha_1+\alpha_1'=(k\alpha_1+r_2) \alpha_2
+\alpha_2'$.

\begin{lemma} \label{edge}
  $s_k>{\beta_3\over \alpha_3}$ and, for sufficiently
  large $k$, there is an edge of the Farey tessellation from
  ${\beta_3\over \alpha_3}$ to $s_k$.
\end{lemma}

\begin{proof}
If $s_k< {\beta_3\over \alpha_3}$, then there is an intermediate
convex torus in $V_3$ parallel to $\bdry V_3$ with infinite slope,
which contradicts the assumptions of Case 2.  Now, if there is no edge
from ${\beta_3\over \alpha_3}$ to $s_k$, then let $s$ be the greatest
rational number in $({\beta_3\over \alpha_3},s_k)$ with an edge to
${\beta_3\over \alpha_3}$, and let $s'$ be the rational number $>s_k$
with an edge to ${\beta_3\over \alpha_3}$ and $s$.  (Such an $s'$
exists provided $k$ is sufficiently large.)  Observe that the
denominator of $s$ is strictly smaller than the denominator of $s_k$
(in absolute value).  Now, there exists a convex torus in $V_3$ with
slope $s$, contradicting the maximality of the twisting number of a
Legendrian fiber in $A$.
\end{proof}

\begin{claim} \label{rel prime}
  The numerator $k(-\alpha_2\beta_1-\alpha_1\beta_2) +
  (1-r_1\beta_1 -\beta_1' -r_2\beta_2 -\beta_2')$ and the denominator
  $k(\alpha_1\alpha_2) + r_1\alpha_1 +\alpha_1'$ of $s_k$ are
  relatively prime.
\end{claim}

\begin{proof}
If not, $\bdry V_3$ has more than two dividing curves, and the
twisting number of the Legendrian regular fiber on $A$ is not maximal.
\end{proof}

\begin{proof}[Proof of Theorem~\ref{infinite GCS}]
According to Lemma~\ref{edge} and Claim~\ref{rel prime}, we know that the determinant of
$(\alpha_3,\beta_3)$ and $(k(\alpha_1\alpha_2) + r_1\alpha_1
+\alpha_1',k(-\alpha_2\beta_1-\alpha_1\beta_2) + (1-r_1\beta_1
-\beta_1' -r_2\beta_2 -\beta_2'))$ is equal to $1$.  In other words,
the determinant of
$(\alpha_1\alpha_2,-\alpha_2\beta_1-\alpha_1\beta_2)$ and
$(r_1\alpha_1 +\alpha_1',1-r_1\beta_1 -\beta_1' -r_2\beta_2
-\beta_2')$ is equal to ${\alpha_1\alpha_2\over \alpha_3}$.  A
straightforward computation gives
$\alpha_1\alpha_2-\alpha_1-\alpha_2={\alpha_1\alpha_2\over \alpha_3}$,
that is, $\sum_{i=1}^3 {1\over \alpha_i}=1$.  This is
precisely the condition for $V$ to be a torus bundle over $S^1$.
\end{proof}

\s\n
{\em Acknowledgements.}  We thank Francis Bonahon for topological help.

\end{document}